\documentclass[11pt]{amsart}
\usepackage{latexsym}
\usepackage{amsmath,amssymb,amsfonts}
\usepackage{graphics}
\usepackage[all]{xy}

\def\Bbb{\mathbb}

\newtheorem{thm}{Theorem}[section]

\newtheorem{lem}[thm]{Lemma}

\newtheorem{question}[thm]{Question}
\newenvironment{xpl}{\mbox{ }\\ {\bf  Example}\mbox{ }}{
\hfill $\diamondsuit$\mbox{}\bigskip}

\def\blacksquare{\hbox to .60em {\vrule width .60em height .60em}}

\begin{document}

\renewcommand{\theequation}{\thesection.\arabic{equation}}

\title[Connected sums with $\Bbb HP^{n}$ or $CaP^{2}$ and the Yamabe invariant]{Connected sums with $\Bbb HP^{n}$ or $CaP^{2}$\\  and the Yamabe invariant}

\author{Chanyoung Sung}
\date{\today}
\address{Department of Mathematics and Institute for Mathematical Sciences \\
Konkuk University\\
         1 Hwayang-dong, Gwangjin-gu, Seoul, KOREA}
\email{cysung@kias.re.kr}
\thanks{This work was supported by the National Research Foundation(NRF) grant funded by the Korea government(MEST). (No. 2010-0016526, 2010-0001194)} 

\keywords{Yamabe invariant, quaternionic projective space, Cayley
plane, surgery}

\subjclass[2000]{53C20, 58E11, 57R65}

\begin{abstract}
Let $M$ be a smooth closed $4k$-manifold whose Yamabe invariant
$Y(M)$ is nonpositive. We show that $$Y(M\sharp\ l\ \Bbb
HP^k\sharp \ m\ \overline{\Bbb HP^k})=Y(M),$$ where  $l,m$ are
nonnegative integers, and $\Bbb HP^k$ is the quaternionic projective
space. When $k=4$, we also have
$$Y(M\sharp\ l\ CaP^2\sharp \ m\ \overline{CaP^2})=Y(M),$$ where $CaP^2$ is the Cayley plane.
\end{abstract}
\maketitle
\setcounter{section}{0} \setcounter{equation}{0}

\section{Introduction}
The Yamabe invariant is an invariant of a smooth closed manifold
defined using the scalar curvature. Let $M$ be a closed smooth
$n$-manifold. By the well-known solution of the Yamabe problem,
each conformal class of a smooth Riemannian metric on $M$ contains
a so-called \emph{Yamabe metric} which has constant scalar
curvature. Moreover, letting $$[g]=\{\varphi g \mid \varphi:M
\rightarrow \Bbb R^+ \  \textrm{is smooth} \}$$ be the conformal
class of a Riemannian metric $g$, a Yamabe metric of $[g]$
actually realizes
$$Y(M,[g]):=\inf_{\tilde{g} \in
[g]} \frac{\int_M s_{\tilde{g}}\ dV_{\tilde{g}}}{(\int_M
dV_{\tilde{g}})^{\frac{n-2}{n}}},$$ where $s_{\tilde{g}}$ and
$dV_{\tilde{g}}$ respectively denote the scalar curvature and the
volume element of $\tilde{g}$. The value $Y(M,[g])$, which is the
value of the scalar curvature of a Yamabe metric with the total
volume $1$ is the \emph{Yamabe constant} of the conformal class.

In a quest of  a ``best'' Yamabe metric or more ambitiously a
``canonical'' metric on $M$, one naturally takes the supremum of
the Yamabe constants over the set of all conformal classes on $M$.
This is possible, because by Aubin's theorem \cite{aubin}, the
Yamabe constant of any conformal class on any $n$-manifold is
always bounded by that of the unit $n$-sphere $S^n(1) \subset \Bbb
R^{n+1}$, which is $n(n-1)(\textrm{Vol}(S^n(1)))^{2/n}$.

The \emph{Yamabe invariant} of $M$, $Y(M)$, is then defined as the
supremum of the Yamabe constants over the set of all conformal
classes on $M$. This supremum is not always attained, but if it is
attained by a metric which is the unique Yamabe metric with total
volume $1$ in its conformal class, then the metric has to be an
Einstein metric.(\cite{ander}) In general, one can hope a singular
or degenerate Einstein metric leading to a kind of a
``geometrization'' from a maximizing sequence of Yamabe metrics.
It is also noteworthy that the Yamabe invariant is a topological
invariant of a closed manifold  depending only on the smooth
structure of the manifold.

The Yamabe invariant of a compact orientable surfaces is just $4\pi \chi(M)$ where
$\chi(M)$ denotes the Euler characteristic of $M$ by
the Gauss-Bonnet theorem. In higher dimensions, it is not an easy task to compute the
Yamabe invariant. Nevertheless recently there have been much progresses in dimension
$3$ and $4$. In dimension $3$, the geometrization by the Ricci flow gives a lot of
answers, and in dimension $4$, the Spin$^{c}$ structure and the Dirac operator are keys for
computing the Yamabe invariant. In particular, LeBrun \cite{lb1,lb3}  showed that if $M$
is a compact K\"ahler surface whose Kodaira dimension is not equal to $- \infty$,  then
$$Y(M)=-4\sqrt{2}\pi\sqrt{(2\chi+3\sigma)(\tilde{M})},$$ where $\sigma$ denotes the
signature and $\tilde{M}$ is the
minimal model of $M$. Now based on this evidence, one can ask if the blowing-up
does not change the Yamabe invariant of a closed orientable $4$-manifold with nonpostive
Yamabe invariant, namely
\begin{question}\label{ques1}
Let $M$ be a closed orientable $4$-manifold with $Y(M)\leq 0$. Is
there an orientation of $M$ such that $Y(M\sharp\ l\
\overline{\Bbb CP^{2}})=Y(M)$ for any integer $l>0$? What about in
higher dimensions?
\end{question}

Further one can also ask whether the analogous statement holds true for the
``quaternionic blow-up'', i.e. a connected sum with the quaternionic projective space
$\overline{\Bbb HP^{k}}$ with reverse orientation, or even a connected sum with the Cayley plane $\overline{CaP^{2}}$ with reverse orientation.
More generally we prove :
\begin{thm}\label{th1}
Let $M$ be a smooth closed $4k$-manifold with $Y(M)\leq 0$.
Then $$Y(M\sharp\ l\ \Bbb HP^{k}\sharp \ m\ \overline{\Bbb HP^{k}})=Y(M),$$ where
$l,m$ are nonnegative integers. When $k=4$, we also have
$$Y(M\sharp\ l\ CaP^{2}\sharp \ m\ \overline{CaP^{2}})=Y(M).$$
\end{thm}

\section{Preliminaries}


A computationally useful formula for the Yamabe constant is
$$|Y(M,[g])|=\inf_{\tilde{g}\in [g]}(\int_{M} {|s_{\tilde{g}}|}^{\frac{n}{2}}
d\mu_{\tilde{g}})^{\frac{2}{n}},$$ where the infimum is attained
only by a Yamabe metric. (For a proof, see \cite{lb3,sung3}.) So
when $Y(M,[g])\leq 0$, this implies that
$$Y(M,[g]) =
-\inf_{\tilde{g}\in [g]}(\int_{M} |s^{-}_{\tilde{g}}|^{\frac{n}{2}}
d\mu_{\tilde{g}})^{\frac{2}{n}},$$ where $s_{g}^-$ is defined as $\min\{s_g,0\}$.
Therefore  when $Y(M)\leq 0$,
\begin{eqnarray}\label{form1}
Y(M)=-\inf_{g}(\int_{M} |s_{g}|^{\frac{n}{2}}
d\mu_{g})^{\frac{2}{n}}=-\inf_{g}(\int_{M} |s_{g}^-|^{\frac{n}{2}}
d\mu_{g})^{\frac{2}{n}}.
\end{eqnarray}
Also essential is Kobayashi's connected sum formula \cite{koba,sung2}
$$
Y(M_{1}\sharp M_{2}) \geq
\left\{
\begin{array}{ll}  -( |Y(M_1)|^{\frac{n}{2}}+  |Y(M_2)|^{\frac{n}{2}} )^{\frac{2}{n}}
   &\mbox{if } Y(M_i)\leq 0 \ \forall i\\
  \min(Y(M_1),Y(M_2)) &\mbox{otherwise,}
\end{array}\right.
$$
which is in fact a special case of the surgery formula in codimension 3 or more (\cite{PY}).
 
We also need to know about the geometry and topology of $\Bbb
HP^{k}$ and $CaP^{2}$. Both have the homogeneous Einstein metric
of positive scalar curvature unique up to constant and can be
viewed as  the mapping cones of the (generalized) Hopf fibrations
$$\pi_{1}: S^{4k-1}\rightarrow \Bbb H P^{k-1}$$ with $S^{3}$ fibers
and $$\pi_{2}: S^{15}\rightarrow S^{8}$$ with $S^{7}$ fibers
respectively. Thus a connected sum with them or their orientation-reversed ones replaces a point with $\Bbb HP^{k-1}$ and $S^8$ respectively so that it deserves the name ``blow-up".

These fibrations have the associated geometries of Riemannian submersion
with totally geodesic fibers.
In case of $\pi_{1}$, $S^{4k-1}$ and $S^{3}$ are endowed with the round metric
of constant curvature $1$, and $\Bbb HP^{k-1}$ is given  the homogeneous Einstein
metric with curvature ranging between $1$ and $4$. In case of $\pi_{2}$, the total
space and the fibers have the round metric of curvature $1$, but the base has
the round metric of  curvature $4$.

We will denote the round $n$-sphere with the metric of constant curvature
$\frac{1}{a^{2}}$ by $S^{n}(a)$, i.e. the sphere of radius $a$ in the Euclidean $\Bbb R^{n+1}$.

\section{Proof of Theorem}
It's enough to prove for one connected sum. Let $M'$ be $M\sharp\
\Bbb HP^{k}$ or $M\sharp \  \overline{\Bbb HP^{k}}$, and set
$n=4k$. First recall that $\Bbb HP^{k}$  admits a metric of
positive scalar curvature meaning that $Y(\Bbb HP^k)>0$. Thus by
the connected sum formula, $$Y(M')\geq Y(M).$$ The idea of the proof
is to surger out an $\Bbb HP^{k-1}$ in $M'$ by performing the
Gromov-Lawson surgery (\cite{GL}) to get back $M$ without decreasing the
Yamabe constant much.

To prove by contradiction, let's assume $$Y(M')>Y(M)+2c > Y(M)$$ for a constant $c>0$.
Let $g$ be an unit-volume Yamabe metric on $M'$ such
that $$s_g\equiv Y(M',[g])= Y(M)+2c.$$ Let $W$ be an $\Bbb
HP^{k-1}\subset \Bbb HP^{k}$ embedded  in $M'$. Take a
$\delta$-tubular neighborhood $$N(\delta)=\{x\in M'|\
dist_{g}(x,W)< \delta\}$$ of $W$ for $\delta>0$. We will take
$\delta$ small enough so that $N(\delta)$ is diffeomorphic to
$\Bbb HP^{k}- \{\textrm{a point}\}$ and the boundary of
$N(\delta)$ is diffeomorphic to $S^{4k-1}$.

First, we consider the case of $Y(M)=0$ so that $s_g>0$.
We perform a Gromov-Lawson surgery described in \cite{GL,sung1,sung2}
on $N(\delta)$ along $W$ keeping the scalar curvature positive to get a cylindrical end isometric to $$(S^{4k-1}\times [0,1],\hat{g}+dt^2),$$ where $(S^{4k-1},\hat{g})$ is a Riemannian
submersion onto $(W,g_W=g|_{W})$ with totally geodesic fibers
isometric to $S^{3}(\varepsilon)$, the round $3$-sphere of radius
$\varepsilon\ll 1$. Here, the horizontal distribution is given by
the connections on the normal bundle. By arranging $\varepsilon$
sufficiently small, $\hat{g}$ has positive scalar curvature.

Now let's take a homotopy
$$H_{b}(t)=\lambda(t)g_W+(1-\lambda(t))g_{std}$$ of smooth metrics
on $W$ from $g_W$ to the homogeneous Einstein metric $g_{std}$ of
$\Bbb HP^{k-1}$ with curvature ranging from $1$ to $4$, where
$\lambda:[0,1]\rightarrow [0,1]$ is a smooth decreasing function
with the property that it is $1$ for $t$ near $0$ and $0$ near
$1$. This induces a homotopy $H_{1}(t)$ of smooth metrics on
$S^{4k-1}$ through a Riemannian submersion with totally geodesic
fibers $S^3(\varepsilon)$. And then we homotope the horizontal
distribution to that of the Hopf fibration through a Riemannian
submersion with totally geodesic fibers $S^3(\varepsilon)$. Let's
denote this homotopy on $S^{4k-1}$ be $H_2(t)$ for $t\in [1,2]$.
When $\varepsilon$ is sufficiently small, $H_{1}(t)+dt^2$ and
$H_{2}(t)+dt^2$ will give a metric of positive scalar curvature on
$S^{4k-1}\times [0,2]$, because it is a Riemannian submersion with
totally geodesic fibers onto $\Bbb HP^{k-1}\times [0,2]$. We
concatenate this part to the above one obtained from the
Gromov-Lawson surgery to get a smooth metric with the boundary
isometric to the squashed sphere $S^{4k-1}$ coming from the Hopf
fibration. Let's denote this metric on the boundary by
$h_\varepsilon$ for a later purpose.

We want to close it up by a $4k$-ball $B^{4k}$ equipped with a
metric of positive scalar curvature. To construct such a metric we
resort to the Gromov-Lawson surgery again. Take a sphere $S^{4k}$
with any metric of positive scalar curvature and let $p$ be any
point on it. As before, we perform a Gromov-Lawson surgery in a
sufficiently small neighborhood of $p$ to get a $4k$-ball with the
positive scalar curvature and the cylindrical end isometric to
$S^{4k-1}(\varepsilon')\times [0,1]$ for a $\varepsilon'>0$. And
then we take a homothety of the whole thing by
$\frac{1}{\varepsilon'}$ so that the boundary is isometric to the
round sphere $(S^{4k-1}(1),h_{1})$. In order to glue this to the
above obtained part, we have to homotope the metric on the
boundary.  We take a homotopy
$$H_{3}(t)=\lambda(t)h_{1}+(1-\lambda(t))h_{\varepsilon}$$ for
$t\in [0,1]$.

\begin{lem}
The metric $H_{3}(t)$ on $S^{4k-1}$ has positive scalar curvature
for every $t\in [0,1]$.
\end{lem}
\begin{proof}
Note that $h_{1}$ and $h_{\varepsilon}$ differ only by the size of
the Hopf fiber. So for each $t$, $H_3(t)$ also has the same
Riemannian submersion structure with the fiber isometric to the
round $3$-sphere of radius
$r(t):=\lambda(t)+(1-\lambda(t))\varepsilon$. By the O'Neill's
formula \cite{besse},
\begin{eqnarray*}
s_{H_{3}(t)}&=&\frac{1}{r^{2}(t)}s_{f}+s_{b}\circ
\pi-r^{2}(t)|A|^{2},
\end{eqnarray*}
where $s_{f}$, $s_{b}$, and $A$ denote the scalar curvature of the
fiber and the base, and the integrability tensor for $t=0$
respectively. Thus $s_{H_{3}(t)}$ is constant for each $t$ and
increases as $t$ increases. From the fact that
$s_{H_{3}(0)}\equiv(4k-1)(4k-2)>0$, the result follows.
\end{proof}

Nevertheless the metric $H_{3}(t)+dt^2$ on $S^{4k-1}\times [0,1]$
may not have positive scalar curvature in general. But due to
Gromov and Lawson's lemma in \cite{GL}, for a sufficiently large
constant $L>0$, $H_{3}(\frac{t}{L})+dt^2$ on $S^{4k-1}\times
[0,L]$ has positive scalar curvature. Now we have a desired
$4k$-ball to be glued to the part made previously out of $M'$.

After the gluing, what we get is just $M$ with a specially devised
smooth metric which we denote by $\bar{g}$. Note that $$s_{\bar{g}}>0,$$ which is a contradiction to the fact that $Y(M)=0$ so that $M$ does not admit a metric of positive scalar curvature metric.

Secondly, in case of $Y(M)<0$, we use a method suggested by C. LeBrun in \cite{lb4}. We first take a
conformal change $e^{2\varphi} g$ of $(M',g)$ such that $\varphi
\equiv 0$ outside $N(\delta)$ and the scalar curvature of
$e^{2\varphi} g$ is positive on a much smaller neighborhood
$N(\delta')$ of $W$.\footnote{One can take $\varphi$ to be a smooth nonpositive function
$\rho(r)$ such that $r(x)=dist_g(x,W)$ for $x\in M'$ and
$$
\rho(r) =
\left\{
\begin{array}{ll}  -ar^2
   &\mbox{for }r\leq \frac{\delta}{3}\\
  0 &\mbox{for } r\geq \delta
\end{array}\right.
$$
where $a>0$ is a  constant. Since $\varphi$ takes the maximum on $W$, $d\varphi$ is identically zero at $W$. For $z\in W$ and a geodesic normal coordinate $(x_1,\cdots,x_n)$ around it such that $\frac{\partial}{\partial x_1}|_z,\cdots,\frac{\partial}{\partial x_{n-4}}|_z$ are tangent to $W$,
\begin{eqnarray*}
(\Delta_g\varphi)|_z &=&(\frac{1}{\sqrt{|g|}}\partial_k(\sqrt{|g|}g^{kl}\partial_l\varphi))|_z= (g^{kl}\partial_k\partial_l(-ar^2))|_z\\ &=&-a(\sum_{i=n-3}^n\partial^2_i(x_{n-3}^2+\cdots +x_n^2))|_z=-8a,
\end{eqnarray*} and hence
\begin{eqnarray*}
s_{e^{2\varphi}g}|_W&=&(e^{-2\varphi}(s_g-2(n-1)\Delta_g\varphi-(n-2)(n-1)|d\varphi|^2))|_W\\ &=& s_g|_W+16(n-1)a,
\end{eqnarray*}
which is positive for sufficiently large $a$.} Moreover one can
arrange that it satisfies
$$-(\int_{M'} |s^{-}_{e^{2\varphi} g}|^{\frac{n}{2}} d\mu_{e^{2\varphi} g})^{\frac{2}{n}} > -(\int_{M'} |s^{-}_{g}|^{\frac{n}{2}} d\mu_{g})^{\frac{2}{n}}-\epsilon $$
for any $\epsilon>0$.\footnote{For this, one replaces the above
$\varphi$ with $\delta^2\rho(\frac{r}{\delta})$ which is in fact equal to $\varphi$ for $r\leq\frac{\delta}{3}$, and takes
$\delta$ sufficiently small. Then as $\delta\rightarrow 0$, the
scalar curvature is bounded below while the volume of $N(\delta)$
tends to zero.} Let's take $\epsilon < c$.

Now we perform a
Gromov-Lawson surgery on $(N(\delta'),e^{2\varphi} g)$ keeping the positivity of the
scalar curvature there. Then in the same way as above, we perform the homotopy process and close it up with a $4k$-ball to get the final metric $\bar{g}$. Since this process makes the scalar curvature positive, we have
\begin{eqnarray*}
-(\int_M |s^-_{\bar{g}}|^{\frac{n}{2}} d\mu_{\bar{g}})^{\frac{2}{n}}
&=&-(\int_{M'} |s^{-}_{e^{2\varphi} g}|^{\frac{n}{2}} d\mu_{e^{2\varphi} g})^{\frac{2}{n}} \\ &>& (Y(M)+2c)-c\\ &>& Y(M).
\end{eqnarray*}
This is a contradiction to the formula (\ref{form1}), and completes a proof for the $\Bbb HP^{k}$ case.

The case of $CaP^{2}$ can be proved in the same way using the fact that $CaP^{2}$ also admits
a metric of positive scalar curvature, and is the mapping cone of the (generalized) Hopf fibration
$\pi: S^{15}\rightarrow S^{8}$ with $S^{7}$ fibers as explained in the previous section.

\section{Examples and Final remarks}

Obviously the theorem is vacuous for the case of $\Bbb HP^{1}$ which is diffeomorphic to $S^{4}$.

\begin{xpl}

Let $H$ be a closed Hadarmard-Cartan manifold, i.e. one with a
metric of nonpositive sectional curvature. By the well-known
theorem of Gromov and Lawson \cite{GL2} on the enlargeable
manifolds, $H$ cannot carry a metric with positive scalar
curvature. Therefore $Y(H)\leq0.$ Applying our theorem to $H$, one
has $$Y(H\sharp\ l\ \Bbb HP^{k}\sharp \ m\ \overline{\Bbb
HP^{k}})=Y(H).$$

For a specific example, take $M=T^{n} \times H$, where $T^{n}$ is
an $n$-dimensional torus and $H$ is as above, e.g. a product of
closed real hyperbolic manifolds. Now since $M$ has an obvious
$F$-structure, its Yamabe invariant is actually $0$ by collapsing the $T^{n}$-part. (Refer to Paternain and Petean \cite{pp}.) Thus  $$Y(M\sharp\ l\ \Bbb
HP^{k}\sharp \ m\ \overline{\Bbb HP^{k}})=0.$$

Similar examples can also be constructed for $CaP^{2}$.
\end{xpl}

Going back to the question \ref{ques1} addressed in the
introduction, our argument does not apply to the case of complex
projective space $\Bbb CP^{k}$. We still have the fact that $\Bbb
CP^{k}$ is the mapping cone of the Hopf fibration $$\pi: S^{2k-1
}\rightarrow \Bbb CP^{k-1}$$ with $S^{1}$ fibers. So the $\Bbb
CP^{k-1}$ is embedded as a submanifold of codimension $2$ which is
one less for the Gromov-Lawson surgery to work. Moreover the
statement corresponding to the theorem \ref{th1} can not be true
at least in dimension $4$. This is because of Wall's stabilization
theorem \cite{wall}. Let $M$ be a simply-connected closed smooth
$4$-manifold.
Then there exists integers  $l,m$ such that $$M \sharp\ l\ \Bbb CP^{2} \sharp\ m\ \overline{\Bbb CP^{2}}=
a\ \Bbb CP^{2} \sharp\ b\ \overline{\Bbb CP^{2}},$$ where $a=l+\frac{1}{2}(b_{2}(M)+ \sigma(M))$
and $b=m+\frac{1}{2}(b_{2}(M)- \sigma(M)).$ But we know that $Y(a\ \Bbb CP^{2} \sharp\ b\ \overline{\Bbb CP^{2}})>0.$
Thus the Yamabe invariant changes drastically by taking connected sums with both $\Bbb CP^{2}$
and $\overline{\Bbb CP^{2}}$. We do not know yet whether the stabilization phenomenon of the Yamabe invariant
is prevalent also in higher dimensions. But at least the question \ref{ques1} is worth investigating in dimension both $4$ and higher.

\end{document}